\documentclass[journal]{IEEEtran}

\usepackage{diagbox}
\usepackage{enumitem}

\usepackage{cite}

\ifCLASSINFOpdf
  \usepackage[pdftex]{graphicx}
\else
\fi
\usepackage{threeparttable}
\usepackage{amsmath}
\usepackage{amsfonts}
\usepackage{bm}

\usepackage{algorithmic}

\usepackage{array}

\ifCLASSOPTIONcompsoc
  \usepackage[caption=false,font=normalsize,labelfont=sf,textfont=sf]{subfig}
\else
  \usepackage[caption=false,font=footnotesize]{subfig}
\fi
\usepackage{fixltx2e}

\usepackage{booktabs}

\usepackage{stfloats}
\usepackage{url}

\usepackage{algorithm,algorithmic}
\usepackage[T1]{fontenc}
\usepackage[latin9]{inputenc}
\usepackage{array}
\usepackage{multirow}
\usepackage{amstext}

\begin{document}

\title{Through-the-Wall Nonlinear SAR Imaging}
%
%
%

\author{Michael~V.~Klibanov,
        Alexey~V.~Smirnov,
        Khoa~Anh~Vo,
        Anders~J.~Sullivan,
        and~Lam~H.~Nguyen~\IEEEmembership{}
\thanks{Michael V. Klibanov is with the Department
of Mathematics and Statistics, University of North Carolina at Charlotte, Charlotte,
NC, 28223 USA (e-mail: mklibanv@uncc.edu).}
\thanks{Alexey V. Smirnov was with the Department
of Mathematics and Statistics, University of North Carolina at Charlotte, Charlotte,
NC, 28223 USA.}
\thanks{Khoa Anh Vo is with the Department
of Mathematics and Statistics, University of North Carolina at Charlotte, Charlotte,
NC, 28223 USA.}
\thanks{Anders J. Sulliavan and Lam H. Nguyen are with he Army Research Laboratory
AMSRD-ARL-SE-RU, Adelphi, MD 20783, USA.}
\thanks{Manuscript received April 19, 2005; revised August 26, 2015. The work was supported by US Army Research Laboratory and US Army
Research Office grant W911NF-19-1-0044.}}


\maketitle

\begin{abstract}
An inverse scattering problem for SAR data in application to through-the-wall imaging is addressed. In contrast with the conventional algorithms of SAR imaging, that work with the linearized mathematical model based on the Born approximation, the fully nonlinear case is considered here. To avoid the local minima problem, the so-called "convexification" globally convergent inversion scheme is applied to approximate the  distribution of the slant range (SR) dielectric constant in the 3-D domain. The benchmark scene of this paper comprises a homogeneous dielectric wall and different dielectric targets hidden behind it. The results comprise two dimensional images of the SR dielectric constant of the scene of interest. Numerical results are obtained by the proposed inversion method for both the computationally simulated and experimental data. Our results show that the values, cross-range sizes and locations of SR dielectric constants for targets hidden behind the wall are close to those of real targets. Numerical comparison with the solution of the linearized inverse scattering problem provided by the Born approximation, commonly used in conventional SAR imaging, shows a significantly better accuracy of our results. 
\end{abstract}


%
\IEEEpeerreviewmaketitle

\section{Introduction}

%
%
%
%

\IEEEPARstart{T}{HROUGH-THE-WALL} imaging (TWI) is an emerging field of technology that addresses the problem of imaging of hidden targets using electromagnetic waves. This problem is of a great interest in a number of civilian missions and has a dual-use with obvious military applications \cite{Nature, Leigs, Lubeck}. 

TWI relies on a stable recovery of target's electromagnetic properties. Among the most valuable for TWI applications is the object's dielectric constant, which is involved as a coefficient of the governing wave-like Partial Differential Equation (PDE). It is well known that the solution of any PDE depends nonlinearly on its coefficient. Hence, the corresponding inverse problem is nonlinear. The complexity of the considered inverse problem gave rise to a variety of imaging algorithms, which can be distinguished by the kind of approximations made. Among the most popular ones are the methods based on the linearization via the Born approximation \cite{Tsynkov,Zhang}, which is valid for weakly scattering targets. It is known that for a general dielectric imaged object the Born-based inversion schemes yield the reconstructions of only locations and shapes. Another type of inversion methods are those directly based the least squares minimization, see, e.g. \cite{Chavent,Gonch}. The common problem with the latter methods, however, is the phenomenon of multiple local minima of corresponding least squares cost functional.

The group of methods, most common in practical TWI applications are migration methods, which are based on the use of the Green's function for solving the inverse problem see, e.g. \cite{Sold, Tiv, Turk} and the references therein. However, any migration algorithm requires an a priori knowledge of the distribution of the dielectric constant to focus and adjust amplitudes of the received signals. Therefore, such an algorithm falls into the category of small perturbation methods, which highly depend on the information about the medium of interest. On the other hand, our research group has demonstrated that our convexification inversion method, mentioned in Abstract, which fully takes into account the inherent nonlinearity of the inverse scattering problem, can achieve a high-resolution quantitative reconstructions of the dielectric properties of targets \cite{Khoa,KlibKol,Klibhyp,SKSN,SKN2}.

Unlike conventional SAR imaging algorithms, which rely on some approximations, e.g. Born approximation, we propose a fully nonlinear reconstruction method, that exploits the same data. Our approach is based on the "convexification" concept, which has been developed by this research group for a number of years, see, e.g. \cite{Kl}-\cite{KT} for initial works and \cite{Khoa,KlibKol,Klibhyp,SKSN,SKN2} and references cited therein for more recent publications. In the convexification, one first changes variables to obtain a boundary value problem for a quasilinear PDE, in which the unknown coefficient is not present. Depending on some specifics, sometimes this might be an integral differential equation or a system of coupled PDEs. The papers on the convexification method for the 1-D coefficient inverse problems for the wave-like PDEs \cite{SKSN}, \cite{SKN2} are particularly close to this one. The key strength of the proposed method is the fact that it is globally convergent, i.e. no a priori knowledge of the solution is needed.

The key step of our method is that we replace the original inverse problem with a number of 1-D inverse problems. Each of these inverse problems is solved via the convexification. We justify our approach by a number of numerical experiments. Another point of the justification comes from the comparison of the raw data after delay-and-sum procedure, used in the preprocessing, with the data, which are computationally simulated for the 1-D wave equation. The solutions of those 1-D inverse problems are then merged to form a 2-D  slant range image of the scene of interest. This image in turn allows us to accurately estimate locations, shapes and permittivities of targets.

The inversion algorithm is further tested on the computationally simulated data. We also demonstrate a good performance of our method on experimental data which were collected for an inspection of a building. We compare the performance of our method with \ the Born approximation method. It is assumed in the Born approximation that $\left\vert \varepsilon _{r}\left( \bm{x}\right) -1\right\vert \ll 1$, where $%
\varepsilon _{r}\left( \bm{x}\right) $ is the dielectric constant. We demonstrate below that the Born approximation method significantly underestimates dielectric constants. Since most current SAR imaging techniques are also based on the Born approximation, then we conjecture that those techniques also significantly underestimate values of dielectric constants of targets, see, e.g. \cite{Amin}, \cite{Sold}.

The remainder of this paper is organized as follows: In section II, we present the mathematical model of TWI and state the inverse scattering problem. In section III we describe the delay-and-sum procedure applied to the raw data and present a novel inversion method for the reconstruction of the dielectric constant of the targets behind the wall. In section IV we provide the results of the numerical tests for both simulated
and experimental data. The performance of our method is compared to that of the linearized inverse scattering problem solved via Born approximation in section V. Finally, we summarize the results of this paper in section VI.

\section{The SAR Data and Our Imaging Goal}

We explain in this section how the SAR data are collected and what kind of image do we want to obtain from these measurements using our inversion procedure. We consider only the so-called "stripmap" imaging configuration \cite{Showman}. One of two main difficulties of SAR imaging is that the data are underdetermined. Indeed, the data depend only on two variables: location of the antenna running along a straight line and time. On the other hand, the unknown dielectric constant $\varepsilon_{r}\left( \bm{x}\right) $ depends on three variables. This is
why a precise mathematical statement of the inverse problem is not possible here. The second main difficulty of the SAR data is their nonlinear dependence on the unknown dielectric constant $\varepsilon _{r}\left( \bm{x}\right) .$

In SAR imaging, a microwave antenna is used to transmit pulses of a certain
duration. A part of the energy is scattered back to the receiver, which, in
general may not be collocated with the antenna, while the remainder is
transmitted into the dielectric medium. Additionally, some energy may be
absorbed by the medium. But we neglect this loss in the present study.
Even
though the complete system of Maxwell's equations governs the propagation of
the EM field in a medium, the conventional mathematical model of SAR for
nonmagnetic, lossless dielectric medium works with a single wave-like PDE (\ref{2.4}) given below \cite{Cheney,Cutrona,Gilman}. 

\begin{figure*}[htb]
\centering
\subfloat[A schematic diagram for transmitter/receiver antenna setup, the elevation angle $\theta $ and the slant range plane $P$.]{\includegraphics[width =.45\textwidth]{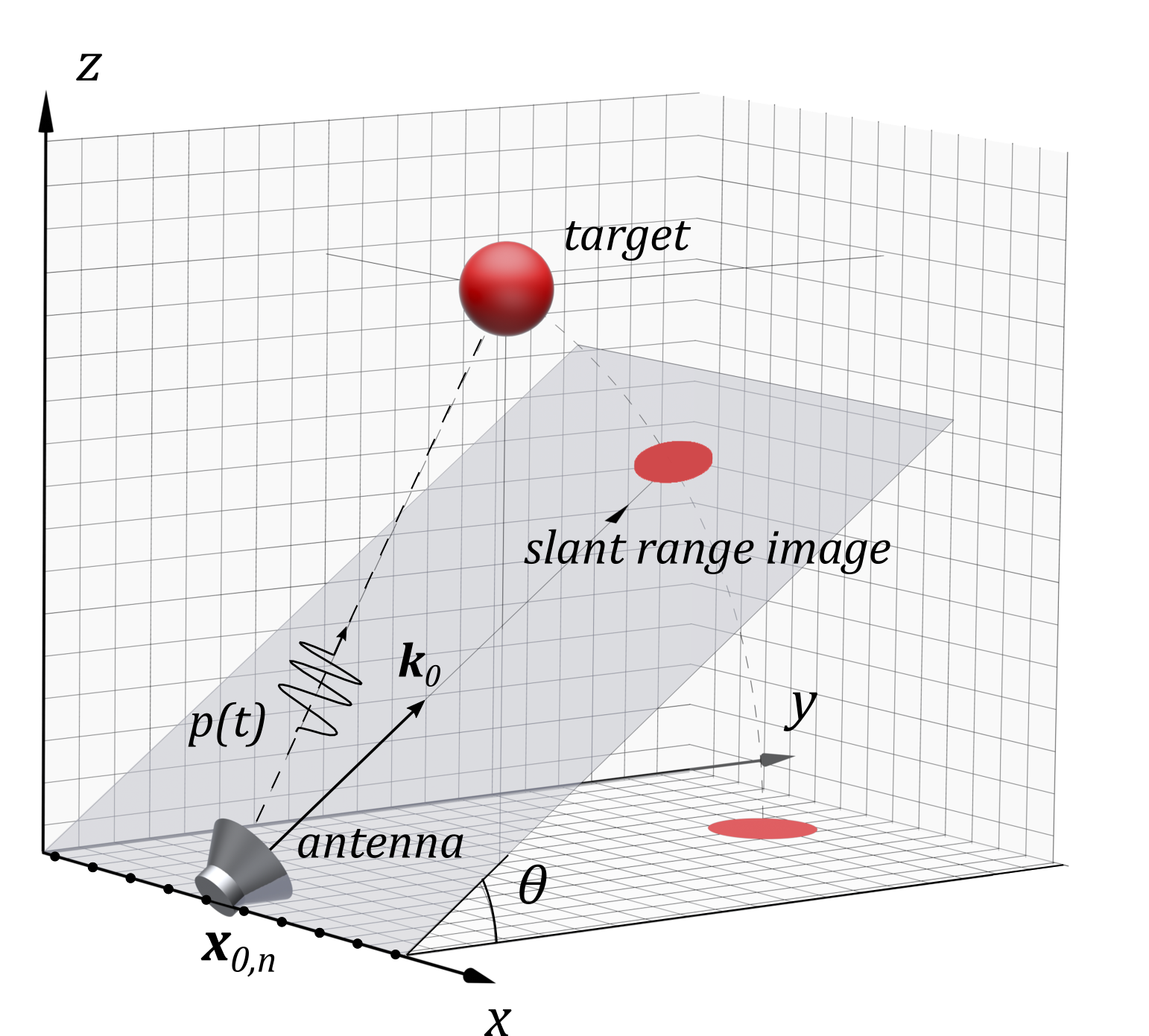}} \hspace{5em} 
\subfloat[Formation of the 2-D
image of the SR distribution of the dielectric constant $\widetilde{\varepsilon}_{r}\left( x, \rho\right)$.]{\includegraphics[width =.36\textwidth]{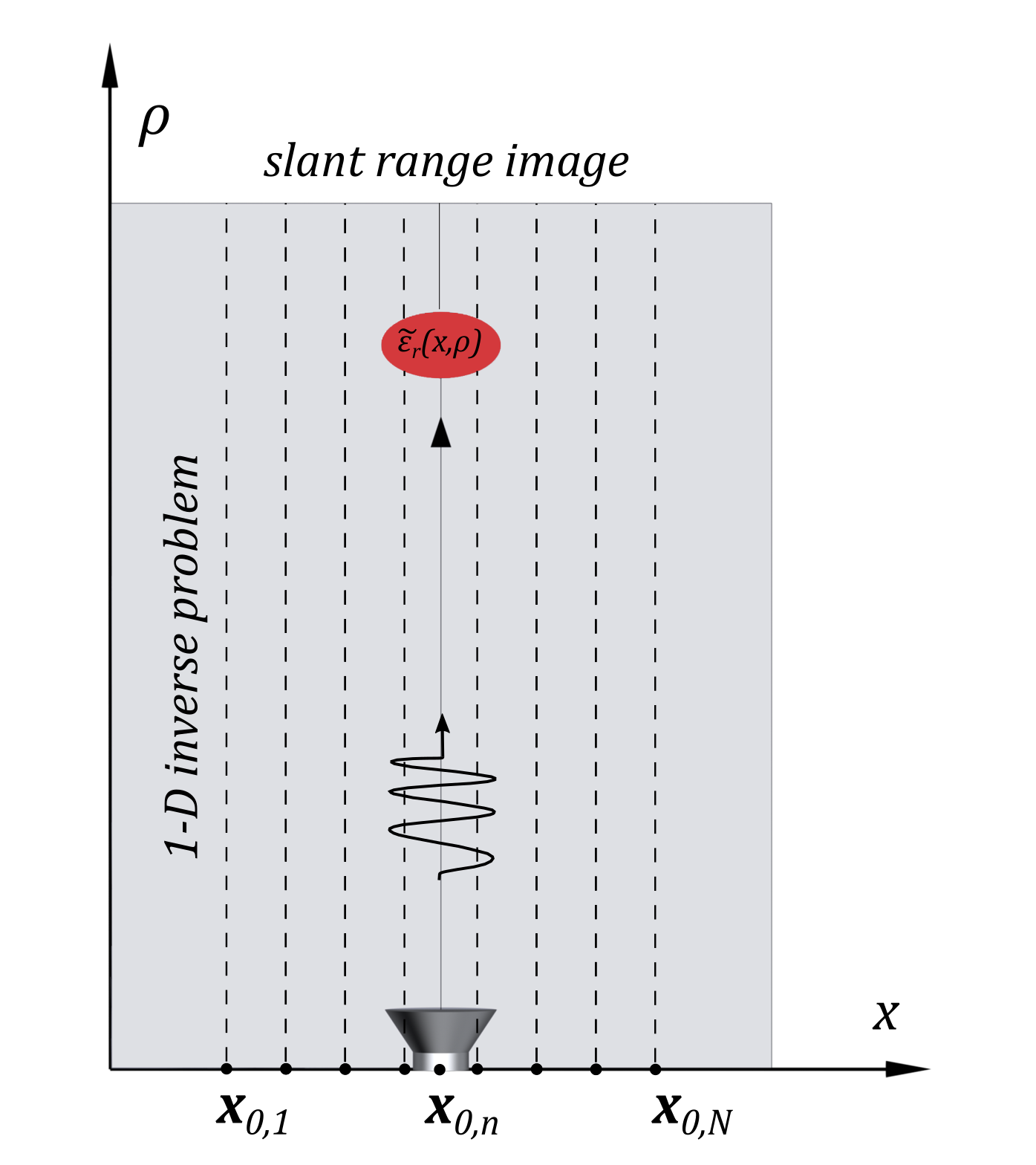}} 
\caption{A schematic diagram of collection of
SAR data in Through-the-Wall imaging problem. Antenna moves along the straight line.}
\label{fig1}
\end{figure*}

\subsection{Measurement Setup}

We provide below definitions of the antennas we use and the slant range. Denote $\bm{x}=\left( x,y,z\right) $ points of the space $\mathbb{R}^{3}.$ Let $\Omega \subset \mathbb{R}^{3}$ be a bounded domain where targets of our interest are located. Below $C^{k},k=0,1,2$ denotes the space of $k$
times continuously differentiable functions. It is assumed that the dielectric constant $\varepsilon
_{r}\left( \bm{x}\right)$ is sufficiently smooth in $\mathbb{R}^{3},$ $\varepsilon _{r}(\bm{x})=1$ outside $\overline{\Omega }$ and $\varepsilon_{r}\left(\bm{x}\right) \geq 1$ inside of $\Omega$. We assume that $\Omega $ is a cube with the center at the point $(x_{c},y_{c},z_{c})$ and the side of length $\widetilde{R}>0.$ We assume that $\Omega $ does not intercept with $x-$axis. Consider  the interval $l_{0}$ of length $L$, along which the transmitting/receiving antenna
runs
\begin{equation*}
l_{0}=\left\{ \bm{x}_{0}=\left( x_{0},0,0\right) \hspace{0.3em} : \hspace{0.3em} x_{0}\in \left(
-L/2,L/2\right) \right\}
\end{equation*}%
This interval is a part of the $x-$axis, chosen to be parallel to the wall, which is
located close to the front face of the domain $\Omega $. Pulses are radiated
at a finite number $N$ of points $\bm{x}_{0,n}\in l_{0},n=1,...,N.$ For each
antenna location, the time resolved backscattering wave is recorded at the
same point $\bm{x}_{0,n}$. Similarly to \cite{Gilman}, we assume that the size of the
transmitting antenna is negligibly small and that the transmitter and the receiver form the same point.
These assumptions works well for large distances between transmitters and the
domain of interest $\Omega .$ At the same time, we use these assumptions only
for the inverse problem. When simulating the data for the forward problem,
we work with a more realistic circular antenna. The speed of the wave
propagation in $\Omega \setminus \mathbb{R}^{3}$ coincides with the speed of
light in the free space $c_{0}$.

We assume that the antenna moves along $x$-axis, we also assume that it is oriented in such a way that the main part of the radiated energy propagates in the direction given by the vector $\bm{k}_0 = (\bm{k}_{0,x},\bm{k}_{0,y},\bm{k}_{0,z})$ and this vector is orthogonal to the $x$-axis. The angle $\theta$ that $\bm{k}_0$ forms with the $x,y-$plane is called the "elevation angle", see Figure 1(a). Consider the plane $P$ which is passing through and $l_0$ and parallel to $\bm{k}_0$. Then $P$ is called the slant range (slant range plane), see Figure 1. The center of the disk of the antenna is located at $\left\{ \bm{x}_{0,n}\right\} $. Denote that disk as $S\left( 
\bm{x}_{0},\theta ,D\right)$ where $D>0$ is the diameter of this disk. Let the number $\eta \in \left( 0,D/2\right) .$ Consider the smoothing function $m\left( \eta ,\theta ,\bm{x},\bm{x}%
_{0},D\right) ,$ which belongs to $C^{2}\left( \mathbb{R}^{3}\right) $ with
respect to $\bm{x}\in \mathbb{R}^{3},$ and is defined as%
\begin{equation}
m\left( \eta ,\theta ,\bm{x},\bm{x}_{0,n},D\right) =\left\{ 
\begin{array}{lcl}
1, \quad \left\vert \bm{x}-\bm{x}_{0,n}\right\vert <D/2-\eta, \\ 
0, \quad \left\vert \bm{x}-\bm{x}_{0,n}\right\vert \geq D/2, \\ 
\in \left[ 0,1\right], \quad \text{otherwise} 
\end{array}
\right.
\label{2.1}
\end{equation}

Hence, $m\left( \eta ,\theta ,\bm{x},\bm{x}_{0,n},D\right) =0$
outside of the ball with the center at $\bm{x}_{0,n}$ and the radius $%
D/2 $.

\subsection{The Forward Problem}

Let $\omega _{0}$ be the carrier (central) frequency of the transmitted signal. For $\tau >0$ define the cut off
function $\chi _{\tau }(t)$ as%
\begin{equation}
\chi _{\tau }(t)=\left\{ 
\begin{array}{lcl}
1,\quad t\in \left( 0,\tau \right), \\
0, \quad \text{otherwise}
\end{array}%
\right.  \label{2.100}
\end{equation}%
Define 
\begin{equation}
p\left( t\right) =\chi _{\tau }(t)e^{-i\alpha \left( t-\tau /2\right)
^{2}}e^{-i\omega _{0}t}  \label{2.2}
\end{equation}%
Expression (\ref{2.2}) for $p\left(
t\right) $ is called "linear modulated pulse" or "chirp"
and $\alpha $ is called the "chirp rate" \cite{Gilman}.

To generate the data for the inverse problem, we work below with the forward problem of finding the funcion $u\left( \bm{x},\bm{x}_{0,n},t\right) 
$ for $\bm{x}\in \mathbb{R}^{3},t\in \left( 0,T\right) $ which satisfies the following conditions:
\begin{eqnarray}
\varepsilon _{r}\left( \bm{x}\right) u_{tt}=\nabla _{\bm{x}%
}^{2}u+Q\left( t,\theta ,\bm{x},\bm{x}_{0}\right), \label{2.3} \\
u\left( \bm{x},\bm{x}_{0,n},0\right) =u_{t}\left( \bm{x},\bm{%
x}_{0,n},0\right) =0,  \label{2.4} \\
Q\left( t,\theta ,\bm{x},\bm{x}_{0,n}\right) =p\left( t\right) 
\widetilde{Q}\left( \theta ,\bm{x},\bm{x}_{0,n}\right)  \label{2.5}
\end{eqnarray}
where $\left( 0,T\right) $ is a certain time interval. The function $\widetilde{Q}\left( \theta ,\bm{x},\bm{x}%
_{0,n}\right) $ is defined below. To define the function $\widetilde{Q}\left( \theta ,\bm{x},\bm{x}%
_{0,n}\right) ,$ we rotate coordinates $\bm{x}=\left( x,y,z\right)
\rightarrow \bm{x}^{\prime }=\left( x^{\prime },y^{\prime },z^{\prime
}\right) $ so that in this new coordinate system the vector $\bm{k}_0$ is aligned with the axis $z^{\prime
} $. In new coordinates $\bm{x}_{0,n}$ becomes $\bm{x}_{0,n}=\left(
x_{0,n}^{\prime },y_{0,n}^{\prime },z_{0,n}^{\prime }\right) \in l_{0}.$ Thus, we define the function $\widetilde{Q}\left( \theta ,\bm{x},\bm{%
x}_{0,n}\right) $ in (\ref{2.3}) as%
\begin{equation}
\widetilde{Q}\left( \theta ,\bm{x},\bm{x}_{0,n}\right) =\delta
\left( z^{\prime }-z_{0}^{\prime }\right) m\left( \eta ,\theta ,\bm{x}-%
\bm{x}_{0,n},D\right) ,  \label{2.6}
\end{equation}%
where $\delta \left( z^{\prime }-z_{0}^{\prime }\right) $ is the delta
function. Hence, $\widetilde{Q}\left( \theta ,\bm{x},\bm{x}%
_{0,n}\right) =0$ outside the disk $S\left( \bm{x}_{0},\theta ,D\right) $%
, i.e. outside of our dish antenna.

To solve the forward problem (\ref{2.1})-(\ref{2.6}) numerically, we apply
to the function $u$ Fourier transform with respect to $t$, obtain the
Helmholtz equation for the resulting function $v$ and then obtain an analog
of the Lippmann-Schwinger equation \cite{Colton},%
\begin{align}
\begin{split}
&v\left( \bm{x},\bm{x}_{0},k\right) =v_{0}\left( \bm{x},\bm{x}%
_{0},k\right) \\ &+k^{2}\int\displaylimits_{\Omega }\frac{\exp \left( ik\left\vert 
\bm{x-\eta }\right\vert \right) }{4\pi \left\vert \bm{x-\eta }%
\right\vert }\left( \varepsilon _{r}\left( \bm{\eta }\right) -1\right)
v\left( \bm{\eta },\bm{x}_{0},k\right) d\bm{\eta }, \label{2.7}
\end{split}
\end{align}
\begin{align}
\begin{split}
&v_{0}\left( \bm{x},\bm{x}_{0},k\right) = X\left( k\right) \times \\
&\times \int\displaylimits_{S\left( \bm{x}_{0},\theta ,D\right) }\frac{\exp \left(
ik\left\vert \bm{x-\eta }\right\vert \right) }{4\pi \left\vert \bm{%
x-\eta }\right\vert }m\left( \eta ,\theta ,\bm{\eta }-\bm{x}%
_{0},D\right) d \bm{\eta}  \label{2.70}
\end{split}
\end{align}
where $X\left( k\right) $ is the Fourier transform of the function $p\left(
t\right) $ in (\ref{2.2}). We solve equation (\ref{2.7}) for $k\in \left[
k_{\min },k_{\max }\right] $ using the numerical method described in \cite{Lecht}. Next, we apply the inverse Fourier transform to the function $%
v\left( \bm{x}_{0n},\bm{x}_{0n},k\right) $. Integration is carried
out over the interval $k\in \left[ k_{\min },k_{\max }\right] .$ This
interval is chosen numerically. As a result, we obtain that SAR data is $\bm{F}(t) = [F_1(t), \dots ,F_N(t)]$,
\begin{equation}
F_n(t) =u\left( \bm{x}_{0,n},\bm{x}%
_{0,n},t\right) ,\quad n=1,\dots,N  \label{2.9}
\end{equation}%
where the function $u\left( \bm{x},\bm{x}_{0,n},t\right) $ is the solution of problem (\ref{2.1})-(\ref{2.6}). 

As it was stated in the beginning of section II, one of the two main difficulties we face is the underdetermined nature of the SAR data (\ref{2.9}): the vector function $\bm{F}(t)$ depends on two variables: one is the discrete variable $n$ and the second one is time $t$. On the other hand, the unknown coefficient $\varepsilon _{r}\left( \bm{x}\right) $ depends on three. This is why we formulate our goal not in a precise mathematical way, as, e.g. in \cite{BK,Khoa,Klibhyp}. We assume that values of $\varepsilon _{r}\left( \bm{x}\right) $ in the domain of interest are unknown, and we assume that $\varepsilon _{r}\left(\bm{x}\right) =1$ both outside of dielectric targets and the wall. The value of $\varepsilon _{r}\left( \bm{x}\right) $ inside the wall is also
unknown. Another assumption we use is that the value of $\varepsilon _{r}\left( \bm{x}\right) $ does not change neither within the wall nor within each target.

The main goal of the present study is to image a certain function $\widetilde{\varepsilon }_{r}\left( \bm{x}\right) $ on the slant range plane $\bm{x}\in P.$ This function $\widetilde{\varepsilon }_{r}\left( \bm{x}\right) $ should characterize well both the value of the dielectric constant and the location of that target, i.e. the distance between the interval $l_0$ and the target. It is also desirable that the cross range size of the support of $\widetilde{\varepsilon}_{r}\left( \bm{x}\right) - 1$ would be close to the real cross range sizes of the target. We call the function $\widetilde{\varepsilon }_{r}\left( \bm{x}\right)$ the slant range (SR) dielectric constant. For each $\bm{x}\in P$ we denote $\bm{x} = (x,\rho)$, where $x$ is the $x$-coordinate of $\bm{x}$ and $\rho$ is the distance between $\bm{x}$ and $l_0$. Therefore SR dielectric constant is denoted below by $\widetilde{\varepsilon}_{r}\left( x, \rho\right)$. Now, suppose that a target does not intersect with the slant range plane $P$. Then it is still possible to image this target, since the signal of our antenna would be reflected back from that target, see Figure 1(a). 

\section{Inversion Method}

In this section, we introduce our inversion method to recover the slant range distribution of the function $\widetilde{\varepsilon} _{r}\left( x, \rho\right) $, mentioned
above. Our approach combines the ideas of conventional through-the-wall imaging using SAR data \eqref{2.9} and globally convergent
numerical method for the inverse scattering problem \cite{SKSN,SKN2}. The algorithm is most easily explained as consisting of three stages that we describe below. Our delay-and-sum procedure is a modification of the delay-and-sum procedure of the conventional SAR imaging \cite{Gilman}.
\vspace{1em}

\begin{itemize}
    \item Stage 1. Apply delay-and-sum to the $N-$D vector of time-dependent data $\bm{F}(t)=\left[ F_{1}(t),F_{2}(t),\dots ,F_{N}(t) \right]$, which is either simulated by the solution of the forward problem (\ref{2.1})-(\ref{2.6}) or collected experimentally. Then we obtain the vector of preprocessed data $\widetilde{\bm{f}}\left( t\right) =[ \widetilde{f}_{1}\left( t\right) ,\dots,\widetilde{f}_{N}\left( t\right) ].$
    
    \item Stage 2. Solve a 1-D ISP with the $\widetilde{f}_{n}\left( t\right)$ as data via a version of \cite{SKSN} of the convexification method for each position of the transmitter/receiver $\bm{x}_{0,n}$. 
    \item Stage 3. Apply the filter of Algorithm 1 to the inversion results $\bm{r}\left( \xi\right) = \left(r_{1}(\xi),r_{2}(\xi),...,r_{N}(\xi)\right)$ at each antenna position to obtain the vector function $\widetilde{\bm{\varepsilon}}_r(\rho)=\left(\widetilde{\varepsilon}_r^1(\rho),\widetilde{\varepsilon}_r^2(\rho),\dots, \widetilde{\varepsilon}_r^N(\rho)\right) $. Then, the resulting coefficient $\widetilde{\varepsilon}^n_r(\rho) = \widetilde{\varepsilon}_r(x_{0,n},\rho)$ for $\rho \in [\rho_{min}, \rho_{max}] $ gives us the 1-D distribution of the desired function at the corresponding transmitter/receiver $\bm{x}_{0,n}$. Finally, the 2-D slant range image of the function $\widetilde{\varepsilon}_r(x,\rho)$ for $\bm{x} \in P$ is compiled from merging of $N$ 1-D distributions $\widetilde{\varepsilon}^n_r(\rho)$.
\end{itemize}
\hspace{1em}

The solutions of above 1-D inverse problems can be obtained independently. Since the assumption that the delay-and-sum procedure provides the data for 1-D inverse problem is not derived and is not based on any physical law, then the proposed approach is a completely heuristic one. At the same time, the global convergence of the convexification method is proved rigorously, see Figure 2.

\subsection{The Delay-and-Sum Procedure}

Let $\Omega $ be the domain of interest as defined above. Denote by $\rho_{min},\rho_{max}$ the minimal and maximal distances between the interval $l_{0}$
and the points in the domain $\Omega $. Suppose that the data $\bm{F}(t)$
are collected for $N$ antenna positions for $t\in \lbrack 0,T]$, where $%
T=2\rho_{max}/c_{0}$. Then the delay-and-sum procedure is as follows: assuming
that the dielectric constant is unit everywhere in $\mathbb{R}^{3}$, a
one-to-one correspondence between the travel-time of the transmitted signal $t$ and the distance $\rho\in [\rho_{min},\rho_{max}]$ is established. The
travel-times can be computed precisely. Then for every point $\bm{x}=(x,y,z) \in P$, for every transmitter position $\bm{x}_{0,n}=(x_{0,n},0,0),$ and the distance between them $\Delta \rho(\bm{x},\bm{x}_{0,n}) =  \sqrt{(x-x_{0,n})^{2}+y^{2}+z^{2}}$ define the corresponding delay time as $$\tau_d(\bm{x},\bm{x}_{0,n}, t) = t \left( \sqrt{1+\left( \frac{2 \Delta \rho(\bm{x},\bm{x}_{0,n})}{c_0 t}\right)^2} - 1\right)$$
Then for any given $t \in [0,T]$ the output of the delay-and-sum procedure is given by 
\begin{equation}
\widetilde{f}_{n}(t)=\frac{1}{N}\sum_{i=1}^{N} I_{i,n}(t) F_{i}(t + \tau_d(\bm{x}_{0,i},\bm{x}_{0,n}, t))
\label{3.1}
\end{equation} 
where $I_{i,n}(t)$ is the indicator
function, showing whether the signal received by the antenna located at  $\bm{x}_{0,n}$ at the moment of time $t$ was picked by the receiver located at $\bm{x}_{0,i}$
\begin{equation}
I_{i,n}(t)=\left\{ 
\begin{array}{lcl}
1,\quad \text{if}\quad \tan ^{-1}{\left( \frac{|\bm{x}_{0,i}-\bm{x}_{0,n}|}{|\bm{x}_{0,n}-\bm{x}_{1,n}|}\right) }<\theta _{0}, &  &  \\ 
0,\quad \text{otherwise} &  & 
\end{array}%
\right.  \label{3.2}
\end{equation}%
where $\bm{x}_{1,n} = \bm{x}_{1,n}(t) = \bm{x}_{0,n} + c_0 t \bm{k_0} / 2$ and $\theta _{0}=1.02\hspace{0.3em}l_{min}/D=2.04\hspace{0.3em}\pi c_{0}/(\omega _{0}D)$ is the half beamwidth of the antenna main lobe \cite{Silver}. The characteristic cone, described by the indicator function $I_{i,n}(t)$, corresponds to the main lobe of the circular dish
antenna. The process, described in \eqref{3.1}-\eqref{3.2}, is performed for
all $N$ antenna positions. The implementation of the procedure
is straightforward. 

In contrast to the conventional SAR imaging algorithms, which perform matched filtering before delay-and-sum \cite{Gilman}, we apply only the delay-and-sum procedure to the raw data. Denote $\bm{f}(t) = Re (\widetilde{\bm{f}}(t))$. Then, given the vector valued function $\bm{f}(t) $, pulse duration $\tau $ and the elevation angle $\theta $ of antenna we first solve the 1-D inverse problems for each antenna position. Even though we propose to use a filtering, but only as a part of the postprocessing applied to the result of the inversion. More precisely, we filter out $N_{0}$ Gauss-like object signal response in the computed function $\bm{r}(\xi) = (r_1(\xi),r_2(\xi),\dots,r_N(\xi))$, (see \eqref{3.140} and subsection C below), assuming that we are supposed to image $N_{0}$ targets behind the wall. The wall is one of those targets. This procedure is sequentially applied to the functions $r_n(\xi)$, computed via convexification at each position $\bm{x}_{0,n}$. See Algorithm 1 for the the detailed description of the filtering. The resulting function $\widetilde{\varepsilon}_r(x,\rho)$ is then used to generate the 2-D image of the scene of interest.

\begin{algorithm}
 \caption{Filter for Gauss-like Objects from \eqref{3.14}-\eqref{3.140}}
 \begin{algorithmic}[1]
 \renewcommand{\algorithmicrequire}{\textbf{Input:}}
 \renewcommand{\algorithmicensure}{\textbf{Output:}}
 \REQUIRE $\xi, r^{comp}_n(\xi)$, $\tau$, $\theta$ $N_0$
 \ENSURE  $\widetilde{\varepsilon}_r^n(\rho)$
  \STATE find $5N_0$ local maximums $\bm{m}$ in $r^{comp}_n(\xi)$
  \STATE save the indices of their positions in $\xi$ to $\bm{p}$ 
  \FOR {$i = 1$ to $5N_0$}
  \FOR {$j = 1$ to $size(max[:,1])$}
  \IF {($ i \neq j$) $\&$ ($\vert max[i,1] - max[j,1] \vert < res$) $\&$ ($max[i,1] > 0$) \& ($max[j,1] > 0$)}
  \STATE $\bm{m}[j] \gets 0$
  \ENDIF
  \ENDFOR
  \ENDFOR
  \STATE $max[N,1] \gets sort(max[N,1])$
  \FOR {$i = 1$ to $N$}
  \STATE 
  \IF {($i \ne 0$)}
  \STATE 
  \ENDIF
  \ENDFOR
 \RETURN $\widetilde{\varepsilon}_n(\rho)$ 
 \end{algorithmic} 
 \label{alg1}
 \end{algorithm}

\subsection{The 1-D versus Delay-and-Sum Data}

In this subsection we provide a
numerical justification of our idea to To form the 2-D image of the SR dielectric constant via merging of $N$ solutions of 1-D ISPs. Temporary  denote $\widetilde{\varepsilon}_r(x_{0,n},\rho) = b(\rho)$. Consider a 1-D analogue of forward problem \eqref{2.1}-\eqref{2.6}
\begin{equation}
\begin{array}{lcl}
b\left( \rho\right) v_{tt}=v_{\rho\rho}+h\left( \rho,t\right) ,\hspace{0.3em} \rho\in \mathbb{R}%
,\hspace{0.3em} t\in \left( 0,T\right) ,  \\ 
v(\rho,0)=v_{t}(\rho,0)=0, \\ 
h\left( \rho,t\right) =\delta (\rho)\hspace{0.3em}e^{-i\alpha \left( t-\tau
/2\right) ^{2}}e^{-i\omega _{0}t}.
\end{array}
\label{3.4}
\end{equation}%
Then, a 1-D analogue of the domain $\Omega $ is the interval $\rho\in (\rho_{\min},\rho_{max})$, where $\rho_{min}>0$, the point $\left\{ \rho=0\right\} $
corresponds to the transmitter position $\bm{x}_{0,n}$ and $\rho_{max}$ corresponds to the furthest distance between $%
\bm{x}_{0,n}$ and $\Omega $. We assume that 
\begin{equation}
b(\rho) = \left\{ 
\begin{array}{lcl}
\in [1,\overline{b}], \quad \rho \in [\rho_{min},\rho_{max}], \\
1, \quad \text{otherwise}
\end{array}%
\right.  \label{3.40}
\end{equation}
where the number $\overline{b}$ is known a priori.

It is the similarity between the solution of \eqref{3.40} and the data $\bm{f}(t)$ after delay-and-sum and filtering which has prompted us to find the function $\widetilde{\varepsilon}_{r}\left( x, \rho \right)$ on the slant range via solving $N$ 1-D ISPs independently and then merging their solutions to obtain the solution to the original ISP, see Figure 1(b). Assume that the transmitted pulse is sufficiently short, i.e. $\tau \ll T$. Then for the formulation of the inverse problem we can approximate the real part $Re\left(h(\rho,t)\right)$ of function in the right hand side of \eqref{3.4} with the delta function as $Re\left(h(x,t)\right) \approx \delta \left( \rho\right) \delta \left( t\right)$. Then the theory of distributions \cite{Vlad} tells us that the initial
conditions in the second line of (\ref{3.4}) should be replaced with $v(\rho,0)=0,v_{t}(\rho,0)=\delta \left( \rho\right) ,$ see the next section.

The numerical solution of \eqref{3.4} for \begin{equation}
    b(\rho) = 1 + 1.5 e^{-4\ln{2}\frac{(x-6.12)^2}{(0.29^2)}} + 4 e^{-4\ln{2}\frac{(x-8.55)^2}{(0.46^2)}}
\label{3.41}
\end{equation}
is compared to the one of \eqref{2.3}-\eqref{2.5} on Figure 2. We took model A of section IV (1-D cross-section of the smoothed dielectric constant of model A, with a straight line, passing through the center of the target) to compare to \eqref{3.41}. Then, the solution $u(\bm{x_{0}},\bm{x_{0}},t), \bm{x_{0}} = (0, 0, 0)$ was multiplied by a calibration factor of $CF = 1.8 \times 10^8$. 

\begin{figure}
    \centering
    \includegraphics[width =.45\textwidth]{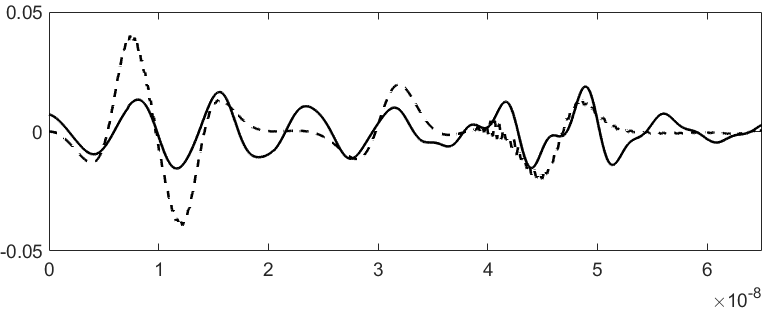}
    \caption{The dashed line shows the solution $v(0,t)$ of the problem \eqref{3.4}. The solid line displays the computationally simulated data $u(\bm{x_{0}},\bm{x_{0}},t)$ for $\bm{x_{0}} = (0, 0, 0)$ after the delay-and-sum preprocessing and scaling (multiplication) with $CF = 1.8 \times 10^{8}$ for model A. The horizontal axis depicts the values of $t$ in seconds.}
    \label{fig2}
\end{figure}

\subsection{1-D Inverse Scattering Problem}

Recall that for each antenna position $\bm{x}_{0,n}$ we use the 1-D variable $\rho \in \left( \rho_{min},\rho_{max}\right) ,$ to describe the distance in the slant range, see the previous subsection. We now scale the interval $ \rho \in (\rho_{min},\rho_{max})$ to the interval $\rho \in (0,1)$ for convenience. Thus, we replace (\ref{3.40}) with%
\begin{equation}
b(\rho) = \left\{ 
\begin{array}{lcl}
\in [1,\overline{b}], \quad \rho \in [0,1] \\
1, \quad \text{otherwise}
\end{array}%
\right.  \label{3.43}
\end{equation}
where the number $\overline{b}$ is the same as in \eqref{3.40}. Then we replace (\ref{3.4}) with: 
\begin{equation}
\begin{array}{lcl}
b\left( \rho\right) U_{tt}=U_{\rho\rho},\hspace{0.3em} \rho\in \mathbb{R},\hspace{0.3em} t\in(0,T), \\ 
U(\rho,0)=0,\quad U_{t}(\rho,0)=\delta (\rho) 
\end{array}
\label{3.44}
\end{equation}%
It was established in \cite{SKSN} that the function $U\left( \rho,t\right) $
satisfies absorbing boundary conditions, 
\begin{equation}
\begin{array}{lcl}
U_{\rho}\left( \rho_{1},t\right) -U_{t}\left( \rho_{1},t\right) = 0,\quad \forall
\rho_{1}\leq 0, \\ 
U_{\rho}\left( \rho_{2},t\right) +U_{t}\left( \rho_{2},t\right) = 0,\quad \forall
\rho_{2}\geq 1.%
\end{array}
\label{3.45}
\end{equation}

1-D ISP. \emph{Let the function} $b\in C^{3}\left( \mathbb{R}%
\right) $\emph{\ satisfies conditions (\ref{3.43}). Determine $b(\rho)$ for }$\rho \in \left( 0,1\right) ,$\emph{\ given the measurements }$f\left( t\right) $\emph{\ and} $g\left( t\right) $ \emph{ for }$t\in \left( 0,T\right)$ 
\begin{equation}
U(0,t)=f\left( t\right) ,\text{ }U_{\rho}(0,t)=g\left( t\right)  \label{3.5}
\end{equation}

Here, for the antenna position number $n$, we denote $%
f\left( t\right) =f_{n}\left( t\right) .$ Note that it is sufficient to know only the function $f\left( t\right) $ since, the first condition (\ref{3.45}) guarantees that 
\begin{equation}
g\left( t\right) =d f(t)/dt = f_t(t)  \label{3.50}
\end{equation}%
To solve this 1-D ISP\ numerically, we implement a version of \cite{SKSN} of
the convexification method. Consider the following change of variables%
\begin{equation}
\xi =\xi \left( \rho\right) =\int \displaylimits_{0}^{\rho}\sqrt{b\left( s\right) }%
ds  \label{3.7}
\end{equation}%
Here, $\xi \left( \rho\right) $ is the travel time which the wave needs to travel to the point located at the distance $\left\{ \rho\right\}$ away from the transmitter/receiver, located at $\left\{ \rho=0\right\} $. The relation (\ref{3.7}) is one-to-one, i.e. 
$\xi =\xi \left( \rho\right) \Leftrightarrow \rho=\rho\left( \xi \right) .$ 

Denote $w\left( \xi ,t\right) =U\left( \rho\left( \xi \right) ,t\right) c^{1/4}\left(
\rho\left( \xi \right) \right) $ and set 
\begin{equation}
S\left( \xi \right) =b^{-1/4}\left( \rho\left( \xi \right) \right) ,\hspace{0.3em} r\left( \xi \right) =\frac{S^{\prime \prime }\left( \xi \right) }{S\left(
\xi \right) }-2\left[ \frac{S^{\prime }\left( \xi \right) }{S\left( \xi
\right) }\right] ^{2}  \label{3.603}
\end{equation}
Hence, the coefficient $r\left( \xi \right) \in C^{1}\left( \mathbb{R}%
\right) $. By (\ref{3.43}) and (\ref{3.603}) 
\begin{equation}
r\left( \xi \right) =0\hspace{0.3em}\text{for}\hspace{0.3em} \xi \in \left\{ \xi <0\right\} 
\hspace{0.3em}\cup \hspace{0.3em}\left\{ \xi >\overline{b}\right\} 
\label{3.604}
\end{equation}%
Using (\ref{3.43}), (\ref{3.44}) and (\ref{3.5})-(\ref{3.603}), we obtain 
\begin{equation}
\begin{array}{lcl}
w_{tt}=w_{\xi \xi }+r\left( \xi \right) w,\hspace{0.3em} \xi \in \mathbb{R}, \hspace{0.3em} t\in [0,T_{1}],  \\
w\left( \xi ,0\right) =0,\quad w_{t}\left( \xi ,0\right) =\delta \left(\xi \right) ,  \\
w\left( 0,t\right) =f\left( t\right) ,\quad w_{\xi }\left( 0,t\right)
=g\left( t\right)  \label{3.9}
\end{array}
\end{equation}
where the number $T_{1}\geq 2\sqrt{\overline{b}}$ depends on $T.$ Thus, it
follows from (\ref{3.7})-(\ref{3.9}) that we have reduced the above 1-D
ISP to the 1-D problem of finding the function $r\left( \xi \right) \in
C^{1}\left( \mathbb{R}\right) $ satisfying (\ref{3.604}).

Let $a\geq \sqrt{\overline{b}}$ be an arbitrary number. Define
the rectangle $R =\left\{ \left( \xi ,t\right) \in \left( 0,a\right) \times
(0,T_{1})\right\} \subset \mathbb{R}^{2}$ and the new function $q(\xi ,t)$ as 
\begin{equation}
 q\left( \xi ,t\right)
=w_{t}\left( \xi , \xi+ t\right)  \label{3.12}
\end{equation}%
Then the second line of (\ref{3.45}) and (\ref{3.7})-(\ref{3.12}) lead to
the non-local boundary value problem for the following nonlinear PDE
\begin{equation}
\begin{array}{lcl}
q_{\xi \xi }-2q_{\xi t}+4q_{\xi }\left( \xi ,0\right) q=0,\hspace{0.3em} \left( \xi
,t\right) \in R, \\
q\left( 0,t\right) =s_{0}\left( t\right) ,\quad q_{\xi }\left( 0,t\right)
=s_{1}\left( t\right) ,\\
q_{\xi }(a,t)=0  \label{3.14}
\end{array}
\end{equation}%
where $s_{0}\left( t\right) =f\left( t\right) ,s_{1}\left( t\right)
=d f(t) / dt+f\left( t\right) .$ We also
have
\begin{equation}
r\left( \xi \right) =4q_{\xi }\left( \xi ,0\right).  \label{3.140}
\end{equation}%
As soon as the function $r\left( \xi \right) $ is found from (\ref{3.140}), the target function $b(\rho)$ can be easily found via calculations, which reverse (\ref{3.7}) and (\ref{3.603}).

\subsection{Convexification for Boundary Value Problem \eqref{3.14}}

Following \cite{SKSN}, we now explain how to obtain an approximate solution of problem (\ref{3.14}) by the convexification method. Consider the function $\varphi _{\lambda }\left( \xi ,t\right) $%
\begin{equation}
\varphi _{\lambda }\left( \xi ,t\right) =e^{-2\lambda \left( \xi +\alpha
t\right) },\quad \alpha \in \left( 0,1/2\right) ,\quad \lambda \geq 1\quad
\label{3.15}
\end{equation}%
where $\alpha $ and $\lambda $ are parameters independent on $\xi ,t$. This
is the Carleman Weight Function for the operator $M=\partial _{\xi
}^{2}-2\partial _{\xi }\partial _{t},$ which is the linear part of the quasilinear Partial Differential Operator (PDO) in the left hand side of the first equation of (\ref{3.14}). In other words, the function $\varphi _{\lambda }\left( \xi ,t\right) $ is involved in the Carleman estimate for the operator $\partial _{\xi
}^{2}-2\partial _{\xi }\partial _{t}$ \cite{SKN2}. Let $R_0>0$ be an arbitrary number. Consider the convex set $Y\left( R_0,s_{0},s_{1}\right) $ 
\begin{equation*}
Y\left( R_0,s_{0},s_{1}\right) =\left\{ 
\begin{array}{lcl}
q\in H^{4}\left( R\right),\hspace{0.3em} \left\Vert q\right\Vert _{H^{4}\left(R\right) }<R_0, \\ q\left( 0,t\right) =s_{0}\left( t\right), \hspace{0.3em}q_{\xi }\left( 0,t\right) =s_{1}\left( t\right) ,\\ q_{\xi }\left( a,t\right) =0,
\end{array}%
\right.
\end{equation*}%
where $H^{4}\left( R\right) $ is a Sobolev space, which is a particular case
of the Hilbert space.

Consider the quasilinear PDO
\begin{equation}
\bm{G}\left( q\right) =q_{\xi \xi }-2q_{\xi t}+4q_{\xi }\left( \xi ,0\right)
q=0\text{ }\left( \xi ,t\right) \in R  \label{3.17}
\end{equation}%
Then, the function $q(\xi ,t)$ can be found from the minimization of the
following weighted Tikhonov-like cost functional for the operator $\bm{G}%
:H^{4}\left( R\right) \rightarrow \mathbb{R}$%
\begin{equation}
J_{\lambda ,\gamma }\left( q\right) =\int \displaylimits_{R}\left[ \bm{G}\left(
q\right) \right] ^{2}\varphi _{\lambda }(\xi ,t)d\xi dt+\gamma \left\Vert
q\right\Vert _{H^{4}\left( R\right) }^{2}  \label{3.18}
\end{equation}%
on the set $\overline{Y\left( R_0,s_{0},s_{1}\right) },$ where bar means that
this is a closed set in the space $H^{4}\left( R\right) $, and $\gamma \in \left( 0,1\right) $ is the regularization parameter. Theorem 4.2 of \cite{SKN2} claims the global strict convexity of functional \eqref{3.18}, since the restrictions on the diameter $2R_0>0$ of the convex set $Y\left(R_0,s_{0},s_{1}\right) $ are not imposed. 

To compute the minimizer of functional (\ref{3.18}), we use the finite difference approximation of the differential operator in (\ref{3.18}) on the rectangular mesh with the step sizes $\bm{\Delta}=(\Delta _{\xi },\Delta _{t})$ and minimize with respect to the vector of values at grid points. Furthermore, we have established numerically in \cite{SKSN} that we can replace in our computations the $H^{4}\left( R\right) -$norm in the penalty term with a discrete analog of a simpler $H^{2}\left( R\right) -$norm and apply the simpler to implement gradient descent method, rather than the conjugate gradient method. Similar approach was used in all our past works about numerical studies of the convexification, see, e.g. \cite{Khoa,KlibKol,Klibhyp,SKSN,SKN2}.

\subsection{The Choices of Parameters}

We have five important parameters to choose: $\lambda ,\gamma ,\alpha
,\Delta _{\xi },\Delta _{t}$.  We have performed  a
cross-validation test w.r.t. $\lambda$ to find its optimal value, all else being constant. The theoretical upper bound for
parameter $\alpha =0.5$ is known from \cite{SKN2}. Thus we set $\alpha=0.49 $ in all further computations. All other parameters were found by trial-and-error. A complete study on the optimal choice of the whole set of them is outside of the scope of this paper. We used models A and A* of section IV as a reference model to scale the data and to simultaneously obtain the values of the parameters, that provide the best possible reconstructed dielectric constant in the target. We found 
\begin{equation*}
\lambda =1.0,\hspace{0.3em}\gamma =10^{-8},\hspace{0.3em}\Delta _{\xi }=0.01,%
\hspace{0.3em}\Delta _{t}=0.02
\label{3.310}
\end{equation*}%
to be optimal values for parameters, when only the position of the wall is known and its dielectric constant is unknown. The data, corresponding to the wall, are not truncated from the measured data.

\begin{equation}
\lambda =2.0,\hspace{0.3em}\gamma =10^{-8},\hspace{0.3em}\Delta _{\xi }=0.01,%
\hspace{0.3em}\Delta _{t}=0.02  \label{3.3}
\end{equation}%
On the other hand, the choice \eqref{3.3} is for the case when the data, corresponding to the wall, are truncated from the measured data. 

We note that even though the theory requires the parameter\textbf{\ }$\lambda $ to be sufficiently large, we found numerically that $\lambda \in [1,3]$ works well. Also, it was established that $\lambda = 2$ is about an optimal number whereas the decrease to $\lambda = 0$ leads to a deterioration of results, see Figure 4 of \cite{SKN2}.

\section{Numerical Results}

In this section we test our inversion algorithm for two measurement scenarios. The measurement setup is depicted in Fig 1(b). All distances mentioned below are measured in meters. 

\subsection{Simulated Data}

In the first case, the length of the
interval $\bm{l}_{0}$ is $L=5.5$ the size of the domain $\Omega $, centered at $(x_{c},y_{c},z_{c})=(0,6.6,3.81)$ is $\widetilde{R}=3.2$. This domain is illuminated by an antenna beam transmitting/receiving chirps from $N=61$ equidistant positions. We use a circular dish antenna of the diameter $%
D=0.7$. The elevation angle $\theta =\pi /6$ and the parameters of the pulse are $\omega _{0}=1885$MHz, $\alpha =1.885\times 10^{17}$, and $\tau = 5$ns. The minimial and maximal wavelengths of the transmitted chirp are $l_{min} = 0.33$ and $l_{max} = 1$ respectively. The wall is a $3.2 \times 0.25 \times 3.2$ dielectric rectangular prism, centered at $(0,5.125,3.81)$, parallel to the $xz$-plane. The dielectric constant  does not change inside of the wall and is equal to $\varepsilon_r$(wall)$ = 2.5$. Hence the distance between the antenna and the wall and the target behind it is at least $5.47$. On the other hand $2D^2/l_{min} = 2.97$, which means that we work in a far field zone \cite{Silver}.

\subsubsection{Reference model}
The data $\bm{F}(t)$ are simulated via the solution of
the problem (\ref{2.1})-(\ref{2.6}). Since we use these data instead of the solution of the 1-D problem \eqref{3.44} and we approximate the function $Re\left(h(x,t)\right)$ with $\delta(\rho) \delta\left(t\right)$ rather than working with the one in (\ref{3.4}),
then we need to scale our data $\bm{f}(t)=\left(f_{1}(t),f_{2}(t),...,f_{N}(t)\right)$ by a calibration factor $CF$. Thus, we propose to use the following calibration model: a wall, with parameters as above and a single ball target behind it (a calibration sphere) see model A below. Next, to find the value of $CF$, we vary it, until the dielectric constant of the target ball, computed via our inversion method, becomes sufficiently close to $2.5$. We use this factor $CF$ in computations for other tested models. Even though the value of CF for the experimental data is different, still a similar
approach is used. More precisely, first, we assume that the dielectric constant of the front drywall, is $\varepsilon_r (wall) = 2.5$ for frequencies ranging from 1GHz to 3GHz, which is a reasonable assumption \cite{Thajudeen}. Next, we find such a number CF that the value of the dielectric constant of the front wall computed by our inversion
procedure becomes sufficiently close to $2.5$. 

We have conducted numerical tests for two types of targets: a ball and a rectangular prism, both hidden behind the wall. The wall is as described above. Let $s_x$ denote the size of the target in cross range (dimension in $x$). Then we describe the geometry of these models below. 

models A and A*. The target is a ball with a diameter of $0.4$ ($s_x = 0.4$), centered at $(0,6.20,4.06)$, the distance between the the line $l_0$ and the center of the ball is $\rho_c = 7.41$. The dielectric constant does not change in the ball and $\varepsilon_r = 2.5$. model A* denotes the same model, but the additional preprocessing was done. More precisely, we have simulated the data for the same domain, but the target ball was absent. This way we can compute the response signal from the wall only. Then these data were subtracted from non-filtered measurements $\bm{F(t)}$ of model A to form the new set of data.

models B and B*. The target is a rectangular prism with sizes $0.9 \times 0.4 \times 0.4$ in $x,y$ and $z$ dimensions accordingly ($s_x = 0.4$). The prism is centered at $(0,6.0,4.61)$, the distance between the the line $l_0$ and the center of the prism is $\rho_c =  7.57$. The dielectric constant does not change within the prism. We have tested two values of the dielectric constant within the prism: $\varepsilon_r=3$ and $\varepsilon_r=4$. For model B* we subtracted the wall clutter from the measurements similarly to model A*.
 
models C and C*. The target is a ball with a diameter of $0.4$, centered at $(0,6.0,4.61)$, the distance between the the line $l_0$ and the center of the prism is $\rho_c = 7.57$. The dielectric constant does not change within the ball. We have tested two values of the dielectric constant: $\varepsilon_r=3$ and $\varepsilon_r=5$. For model C* we subtracted the wall clutter from the measurements similarly to model A*.

\subsubsection{Image Postprocessing and Artifact Removal}

We have noticed that the images of SR the dielectric constant, obtained by our inversion method contain two main groups of artifacts:

\begin{enumerate}[label=\alph*)]
\item Value artifact: the computed value of the $\widetilde{\varepsilon}_r(x,\rho)$ changes significantly inside the imaged target.
\item Shape artifact: the estimated size of the target is smaller than the theoretical resolution of the imaging system.
\end{enumerate}

Assume that the dielectric constant $\varepsilon_r$ is homogeneous inside the imaged target. For a selected region of the slant range image (see bottom images on Figure 4(a),(b)), denote the median of the computed SR dielectric constant by $\widetilde{\varepsilon}_r^{comp}$. Then one can eliminate the first type of artifacts by truncation of 1-D dielectric constant distributions $\widetilde{\varepsilon}^n_r$ for a number of antenna positions, if at each position $\bm{x}_{0,n}$ the relative deviation $\vert \widetilde{\varepsilon}_r^n - \widetilde{\varepsilon}_r^{comp} \vert / \widetilde{\varepsilon}_r^{comp} > \sigma$, where $\sigma = 0.15$ is the deviation from the median value, chosen numerically for the reference model and is fixed for all subsequent numerical tests. We compare the median value with the true value of the dielectric constant $\varepsilon_r$ in the target for all models. The corresponding relative errors (in percent) are denoted by $\delta \widetilde{\varepsilon}_r^{comp}/ \varepsilon_r$ and are given in Tables I and III. The second type of artifacts were tackled by comparing the size of selected regions of image with the theoretical cross range $D/2 = 0.35m$ resolution. The targets of smaller sizes were truncated. The results of the inversion prior to the truncation and artifact removal of this section are compared to the processed images on Figure 4. The processed images were linearly interpolated to a four (4) times finer mesh.  

\subsubsection{Reconstruction Results}

The finite-difference analogue of the functional $J_{\lambda ,\gamma }(q)$ is minimized via the gradient decent method. This is quite computationally intensive for reasonably-sized image sampling. We choose $\Delta_\xi = 0.01, \Delta_t = 0.02$, which allows us to reduce the computational cost for the inversion, while maintaining sufficient resolution. See \cite{SKSN, SKN2} for more details. Figure 3 demonstrates the 2-D cross section of the imaged domain with a plane $z = 4.61$, parallel to the $xy$-plane and passing through the center of the target. Note that the distances are different from the ones of Figure 4, the distances on Figure 3 are measured from the $xz$-plane. since Figure 4 displays images obtained for models B and C. We have not included the images of the wall in these figures because the reconstructed value dielectric constant of the wall is too large. The accurate reconstruction of the wall's dielectric properties should be considered separately. So, we focus only on imaging of unknown targets behind the wall. The following parameters are estimated from the images computed for models B, C: the relative error of the target's size in cross range (target's dimension in $x$) $\delta s_x^{comp} / s_x = \vert s_x^{comp} - s_x \vert / s_x $, where $s_x^{comp}$ is the length in meters, estimated from the image and $s_x$ is the true length in $x$. $\delta \rho_c^{comp} / \rho_c$ is the relative error of the distance between the center of the imaged target and $l_0$, $\rho_c^{comp}, \rho_c$ denote the computed distance and the true distance correspondingly. 

\begin{threeparttable}[htb]
    \caption{Dielectric constants computed via Convexification}
\label{table1}
    \small
    \setlength\tabcolsep{0pt}
\begin{tabular*}{\linewidth}{@{\extracolsep{\fill}} l c c c c c c}
    \toprule
 model & $s_x$ & $\delta s_x^{comp} / s_x$ & $\delta \rho_c^{comp} / \rho_c$ & $\widetilde{\varepsilon}_r^{comp}$ & $\delta \widetilde{\varepsilon}_r^{comp} / \varepsilon_r$  \\
        \midrule
B ($\varepsilon_r = 3.0$)    &0.71    & 21\%  & 4.5\%   & 2.89  & 3.7\%\\
B ($\varepsilon_r = 4.0$)     &0.74   & 17\%  & 2.1\%   & 4.18  & 4.4\%\\
B* (no wall)    &0.74    & 17\%  & 2.1\%  & 3.88  & 3.0\%\\
C ($\varepsilon_r = 3.0$)  &0.36      & 10\% & 5.3\% & 2.84  & 5.3\%\\
C ($\varepsilon_r = 5.0$)    &0.34     & 15\% & 4.9\% & 4.34  & 13.2\%\\
C* (no wall)    &0.34    & 15\% & 4.9\%  & 4.46  & 11.0\%\\
        \bottomrule
        \vspace{1em}
\end{tabular*}
\end{threeparttable}

The relative error of the computed dielectric constant $\delta \widetilde{\varepsilon}_r^{comp} / \varepsilon_r$ is defined similarly to $\delta s_x^{comp} / s_x$. The results for the model A are not present in Table I, since we used it as a reference model for the data scaling.

\begin{figure*}[htb]
\centering
\subfloat[Dielectric constant distribution for model B. Cross-section by the plane $z = 4.61$, through the center of the target.]{\includegraphics[width =.45\textwidth]{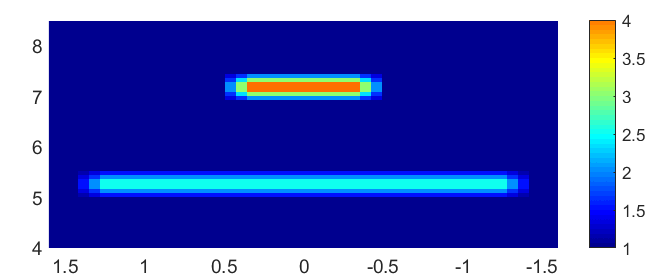}} \quad
\subfloat[Dielectric constant distribution for model C. Cross-section by the plane $z = 4.61$, through the center of the target.]{\includegraphics[width
=.45\textwidth]{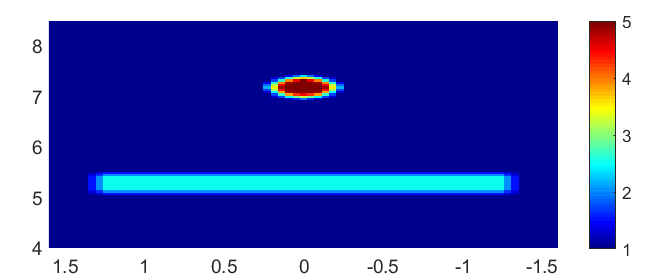}} 
\caption{Images of the the 2-D cross sections of dielectric constants for models B and C. Blue line represents the wall.}
\end{figure*}

\begin{figure*}[htb]
\centering
\subfloat[Results for model B before (top) and after (bottom) postprocessing. ]{\includegraphics[width =.45\textwidth]{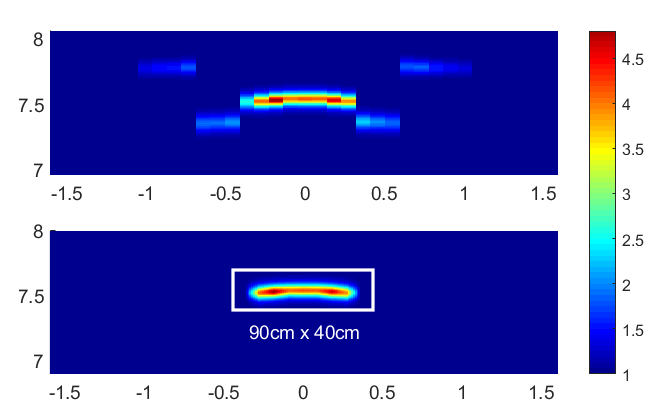}} \quad
\subfloat[Results for model C before (top) and after (bottom) postprocessing.]{\includegraphics[width
=.45\textwidth]{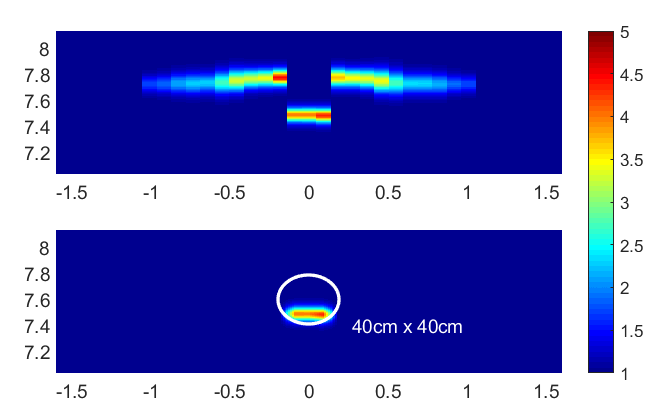}} 
\caption{Slant range dielectric constants for model B with $\varepsilon_r=4$ and for model C with $\varepsilon_r=5$, see section II. The processed images are compared with slant range projections of targets (white line represents the expecetd image of the target, the number shows the true length of the target in cross range). The color scale shows the value of the slant range dielectric constant. The wall clutter is not shown on these images.}
\end{figure*}

\subsection{Experimental Data}

\subsubsection{Measurement Setup}

We use the radar data, experimentally collected by Sullivan and Nguyen with the goal to inspect a building. The tranmitting/receiving setup was driven on the top of a vehicle along two sides of an H-shaped building at a distance of 10 meters parallel to the walls. The antenna takes measurements at $N=2000$ equidistant positions along path 1 and $N= 2000$ positions along path 2, see Figure 5. The parameters of the pulse are $\omega _{0}=1600$MHz, $\alpha = 8.0\times 10^{17}$ and $\tau =1$ns, see \cite{Lam} for a detailed description of the measurement setup.

\begin{figure}
    \centering
    \includegraphics[width =.45\textwidth]{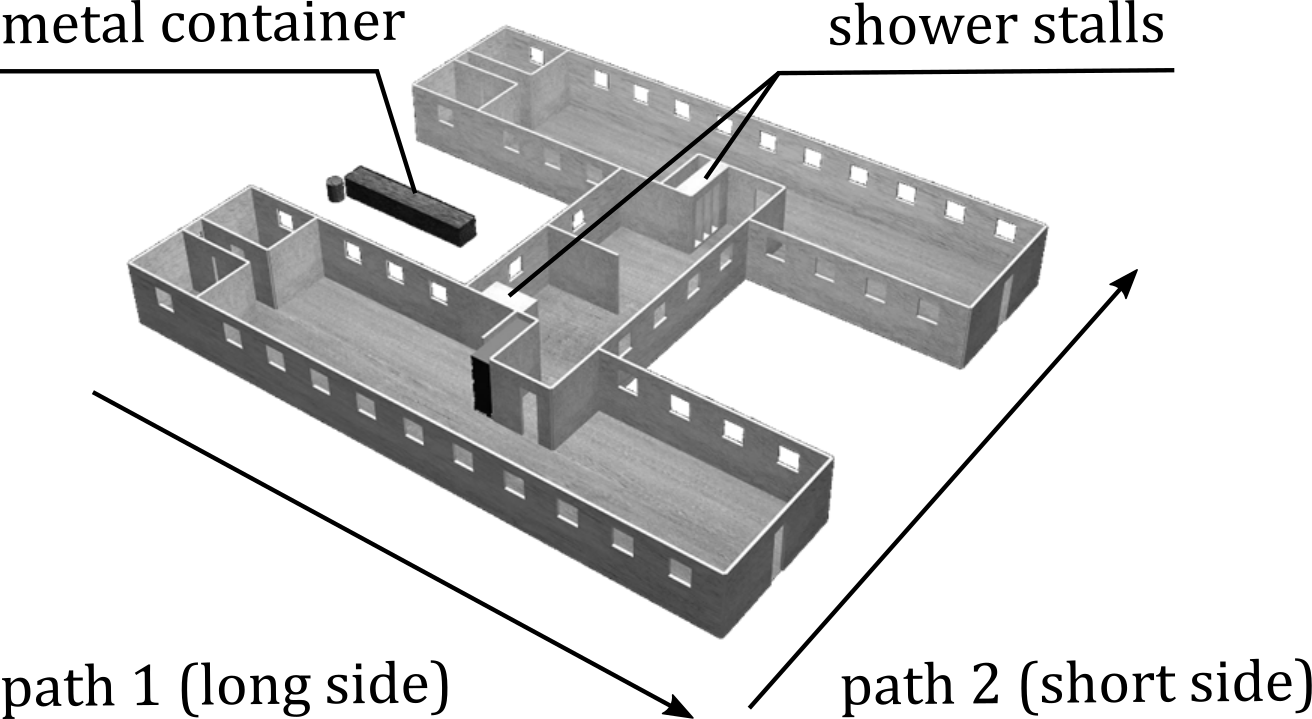}
    \caption{A model of the H-shape barrack building. Two main buildings connected by a short hall way. There are three main imaged objects inside the building: a metallic containerbetween the main buildings and two shower stalls in the hall area.}
    \label{fig5}
\end{figure}

\subsubsection{Results}

As described above, in the case of experimental data, we have chosen the calibration factor CF in such a way that the dielectric constant of the front wall became $\varepsilon_r (wall) = 2.5$ We also note that the conventional SAR image we got does not show the values of the dielectric constants, but rather shows the relative refectivity of the imaged objects. Thus we present the computed contrasts $\nu = \varepsilon_r^{comp} / 2.5$ rather than the computed median values of the SR dielectric constants in the targets of interest, placed inside the building. We compare these contrasts to the ones of the SAR image $\nu_{SAR}$, between the targets and the front wall. Note that these contrasts can further be used to accurately classify types of imaged objects. The inspected building is depicted on the Figure 5. It consists of several drywalls, with ceramic shower stalls and metallic containerhidden behind them. The contrast of the metallic containerin the left had side of Figure 6(a) is about $\nu \approx 7.2$. At the same time, the contrast of the same object on the SAR image of Figure 6(b) is close to $\nu_{SAR} \approx 1.1$. This means that our method allows one to identify the difference between the material of the metallic containerand the front wall, while the SAR image does not. Since the computed $\nu \approx 7.2$ in the metal container, then the computed dielectric constant of this container is $\varepsilon_r^{comp} (container) \approx 7.2 \times 2.5 = 18$. This value can be then used to classify the target. The scale on Figure 6 is logarithmic. The results are summarized in Table II.

Our method assumes that the targets are surrounded by free space. However, we have here the side wall surrounded by the ground, the back wall and two side walls, all of which forming 90 degree angles with each other. This leads to multiple refection waves, all of which come back to the radar, in addition to the direct response from the target wall. Thus, the ground and the walls effectively form a 4-sided corner reflector, which is a quite complicated structure. This is the reason why our technique does not provide good reconstruction results for the side wall, see Figure 6(a). Similar arguments are applicable to the SAR image of Figure 6(b). Another difficulty is to describe the material of the target with a single median value of the computed SR dielectric constant when the target is not homogeneous. This is the case for the shower stalls, see Figures 5 and 6(a). However, we still use the median value of the SR dielectric constant in stalls for simplicity, see Table II, and leave this question for future research.  
\begin{figure*}[htb]
    \centering
    \subfloat[SR dielectric constant of the inspected building, obtained from combination of images along paths 1 and 2. The color scale shows the value of the slant range dielectric constant.]{\includegraphics[width =.43\textwidth]{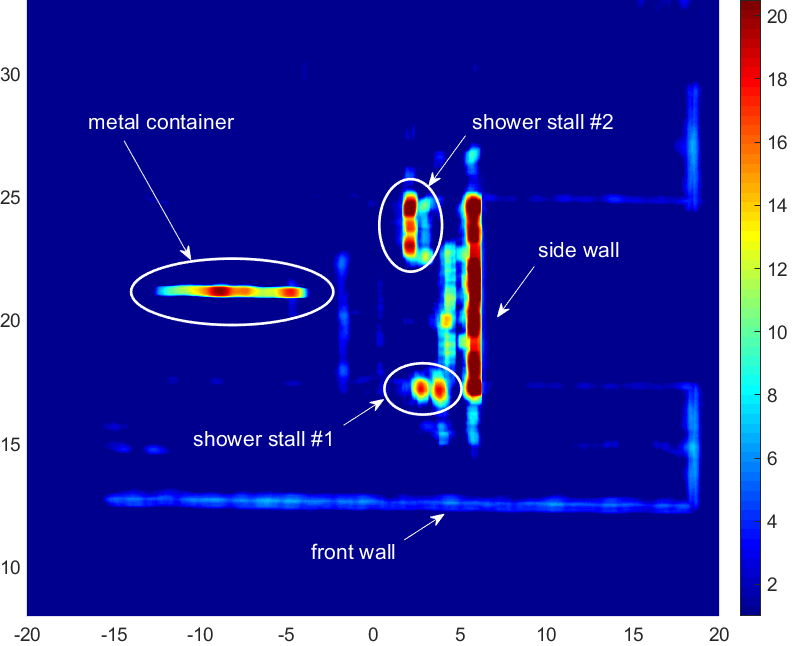}} \quad
    \subfloat[SAR image of the building, obtained from combination of images along paths 1 and 2.]{\includegraphics[width =.43\textwidth]{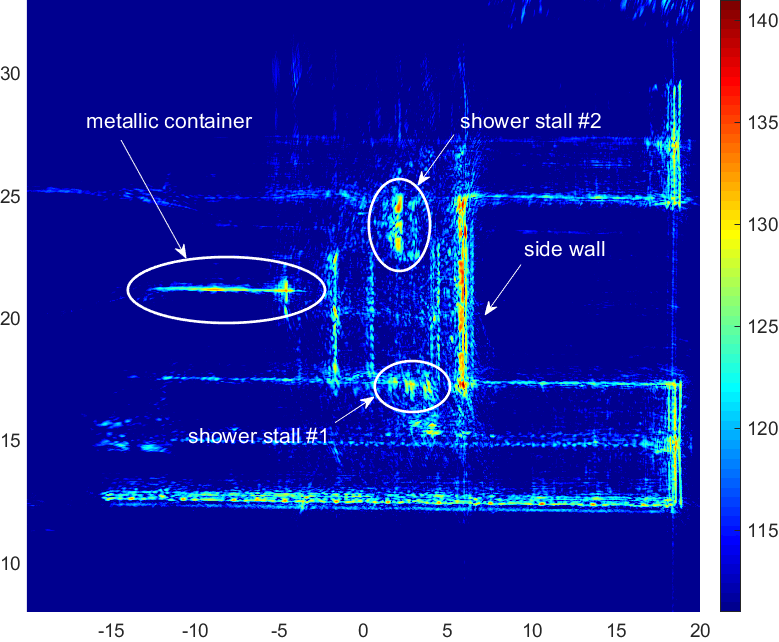}} 
    \caption{Images of the inspected building obtained via our inversion method and SAR.}
\end{figure*}

\begin{threeparttable}[htp]
    \caption{Image Contrast Enhancement for Experimental Data}
\label{table2}
    \small
    \setlength\tabcolsep{0pt}
\begin{tabular*}{\linewidth}{@{\extracolsep{\fill}} l c c c c}
    \toprule
Target & $\widetilde{\varepsilon}_r^{comp}$ & $\varepsilon_{r,table}$ & $\nu$ & $\nu_{SAR}$ \\
        \midrule
Shower stall \#1    & 11.8 & [4,15]\tnote{a} & 4.7 & 1.05 \\
Shower stall \#2    & 14.5 & [4,15]\tnote{a} & 5.8 & 1.08\\
metallic container    & 17.9 & [10,30]\tnote{b} & 7.2 & 1.1 \\
        \bottomrule
\end{tabular*}
\begin{tablenotes}
\item[a] Shower stalls were made of Phenol Formaldehyde Resin and/or Acrylic Resin, the dielectric constant varies greatly, depending on the of type of admixture used, see \cite{Clipper}. \item[b] The so-called "apparent dielectric constant" of a metal is in the interval [10,30], see \cite{Kuzh}.
\end{tablenotes}
\vspace{1em}
\end{threeparttable}
$\varepsilon_{r,table}$ in Table II denotes the interval of values of the dielectric constant for the material of the target, since the true values of the dielectric constants in different materials are given by the interval of values rather than a single value, see \cite{Clipper}. 

\section{Comparison with the Born approximation}

An interesting question is about a comparison of the performance of our nonlinear method with the conventional methods for SAR imaging. All those methods rely on the Born approximation. Furthermore, any of these methods imposes some additional assumptions on the model of Born approximation. Hence, we have decided to compare here our inversion method with the straightforward implementation of Born approximation model. We use a richer set of data in our model than the SAR data, used for inversion in section IV. In other words, the data in this section depends on three rather than two variables. Furthermore, the SAR data \eqref{2.9} are a part of the data used in this section. Thus, it seems to be that the image provided by the Born approximation method of this section should be at least not worse than the one provided by any other version of the SAR imaging algorithm.

It follows from (\ref{2.7}) and (\ref{2.70}) that the inverse scattering problem in the Born approximation case can be stated as: find the function $b\left( 
\boldsymbol{\eta }\right) =\varepsilon _{r}\left( \boldsymbol{\eta }\right)
-1$ from the equation%
\begin{equation}
\int\displaylimits_{\Omega }\frac{e^{ik\left\vert \bm{x-\eta }\right\vert} }{4\pi \left\vert \bm{x-\eta }\right\vert }%
v_{0}\left( \bm{\eta },\bm{x}_{0},k\right) b\left( \boldsymbol{\eta }%
\right) d\bm{\eta =}h\left( \bm{x},\bm{x}_{0},k\right)
\label{5.3}
\end{equation}
where the function $v_{0}\left( \bm{x},\bm{x}_{0},k\right) $ is
defined in (\ref{2.70}), and the function $h\left( \bm{x},\bm{x}%
_{0},k\right) =v\left( \bm{x},\bm{x}_{0},k\right) -v_{0}\left( 
\bm{x},\bm{x}_{0},k\right) .$ We assume that the function $h\left( 
\bm{x},\bm{x}_{0},k\right) $ is known for $\bm{x},\bm{x}%
_{0}\in l_{0}$ and $k\in \left( k_{\min },k_{\max }\right) .$ Since points $%
\bm{x}$ and $\bm{x}_{0}$ run independently over $l_{0},$ then the
data $h\left( \bm{x},\bm{x}_{0},k\right) $ in (\ref{5.3}) depend on
three variables. Furthermore, the $h\left( \bm{x}_{0},\bm{x}
_{0},k\right) $ part of the data $h\left( \bm{x},\bm{x}_{0},k\right) $
actually has the same information content as the the Fourier transform of
the regular SAR data (\ref{2.9}). Consider the set $S=l_{0}\times
l_{0}\times \left( k_{\min },k_{\max }\right) .$ Denote the operator in the
left hand side of (\ref{5.3}) by $A:H^{1}\left( \Omega \right) \rightarrow
L_{2}\left( S\right) $. Assuming that $h\in L_{2}\left( S\right) ,$ Equation
(\ref{5.3}) can be rewritten as 
\begin{equation}
A\left( b\right) =h  \label{5.4}
\end{equation}
It follows from (\ref{5.3}) that (\ref{5.4}) is an integral equation of the
first kind. The latter means that the problem of the solution of (\ref{5.4})
is ill-posed. Hence, we solve it via the regularization method. We minimize
the following functional 
\begin{equation}
J_{\beta }\left( b\right) =\left\Vert A\left( b\right) -h\right\Vert
_{L_{2}\left( S\right) }^{2}+\beta \left\Vert b\right\Vert _{H^{1}\left(
\Omega \right) }^{2} \label{5.6}
\end{equation}%
where $\beta =10^{-4}$ is the regularization parameter. We have found $\beta 
$ numerically by trial-and-error, using model A of section IV as a reference. In doing so, we found such a number $\beta = 10^{-4}$, that maximum of the function $b(\bm{\eta}) = 2.5$ as in model A. The data for the signal reflected from the wall were not counted. Then we have used the same $\beta =10^{-4}$ for model B with a different values of the dielectric constant in the target. The response from the wall, was subtracted from the data $v\left(\bm{x},\bm{x}_{0},k\right) $ for $\left( \bm{x},\bm{x}_{0},k\right) \in S$. The target's median SR dielectric constant $\widetilde{\varepsilon}_r^{comp}$ was computed and compared to both the true one and the value, obtained via our inversion algorithm. Although the shape of the target was reconstructed well, the value of the dielectric constant in it was quite inaccurate one, see Table III.

\begin{threeparttable}[htb]
    \caption{Relative error of the dielectric constant reconstruction via Born Approximation and Our Method}
\label{table3}
    \small
    \setlength\tabcolsep{0pt}
\begin{tabular*}{\linewidth}{@{\extracolsep{\fill}} l cc cc cc @{}}
    \toprule
  &  \multicolumn{2}{c}{model B ($\varepsilon_r = 4.0$)} & \multicolumn{2}{c}{model B ($\varepsilon_r = 6.0$)} \\
        \cmidrule{2-5}
   Method
   & $\widetilde{\varepsilon}_r^{comp}$  
                            & $\delta\widetilde{\varepsilon}_r^{comp} / \varepsilon_r$
                                    & $\widetilde{\varepsilon}_r^{comp}$ 
                                            & $\delta\widetilde{\varepsilon}_r^{comp} / \varepsilon_r$\\
        \midrule
Born approximation          & 2.32  & 42.0\%  & 2.31  & 64.5\%\\
Our inversion       & 4.18  & 4.4\%  & 5.56  & 7.4\%\\
        \bottomrule
        \vspace{1em}
\end{tabular*}
\end{threeparttable}

Thus, the Born approximation leads to much more significant errors in the
value of the dielectric constant than our method. Since SAR imaging algorithms rely on the Born approximation, then we conjecture that the conventional algorithms also lead to high errors in the value of the dielectric constant. 

\section{Conclusion and Discussion}

We have presented a novel nonlinear inversion method for SAR data and
applied it for the reconstruction of the dielectric constant from the
through-the-wall imaging radar data. The proposed method is based on the
numerical solution of a rigorously formulated ISP. A numerical method is
constructed to rigorously ensure the global convergence of the minimization
process. The fact that properties of front walls may be unknown in practice
presents a challenge in producing accurate and reliable images. Our results
for both computationally simulated and experimentally collected data results
confirm that our method can do both: accurately localize targets of
different shapes and accurately compute their dielectric constants. This is true even in the case when the thickness and the dielectric constant of the wall are unknown. The comparison with the Born approximation case shows that our method is significantly more accurate in
computations of dielectric constants. Since any conventional SAR imaging
technique is based on the Born approximation and also uses far less data
than we did in (\ref{5.3})-(\ref{5.6}), then we conjecture that our
technique significantly outperforms any conventional SAR algorithm in the
accuracy of computed dielectric constants of targets in the problem of
through-the-wall imaging.

Thus, the use of the proposed inversion method enhances the capacity of
through-the-wall radar imaging in obtaining images of targets hidden behind
the wall and in their classification.


%

\appendices


\section*{Acknowledgment}

The authors would like to thank Dr. Loc Nguyen and Dr. Mikhail Gilman for many fruitful discussions.

\ifCLASSOPTIONcaptionsoff
\newpage \fi



%

\begin{IEEEbiography}    [{\includegraphics[width=1in,height=1.25in,clip,keepaspectratio]{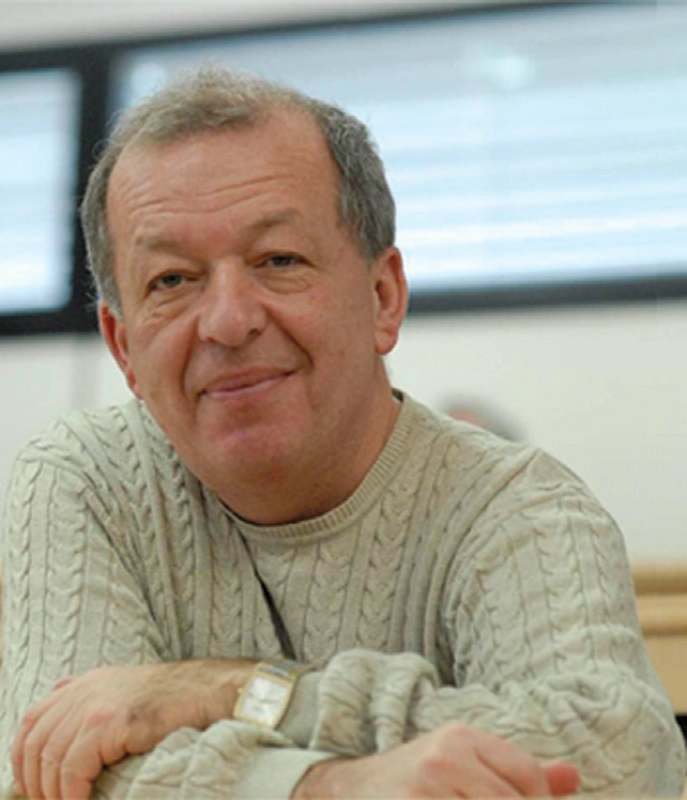}}]{Michael V. Klibanov} has graduated from Novosibirsk State University (NSU), Novosibirsk, Russia, in 1972. NSU is one of very top Russian universities. He got MS in Mathematics. In 1977 he got PhD in Mathematics from Urals State University, Yekaterinburg, Russia. In 1986 he got the highest scientific degree, Doctor of Science in Mathematics from Computing Center of the Siberian Branch of Russian Academy of Science, Novosibirsk. Through his entire career Klibanov works solely on inverse problems. The paper of A.L. Bukhgeim and M.V. Klibanov, "Uniqueness in the large of a class of multidimensional inverse problems" Soviet Mathematic, Doklady, 17, 244-247, 1981 became one of very few classical papers in the field of Inverse Problems. In this paper the powerful tool of Carleman estimates was introduced in the field for the first time. While previously Klibanov has worked only on the uniqueness issue, currently he develops globally convergent numerical methods for Coefficient Inverse Problems without overdetermination. He has published a total of 170 papers and his works were cited 2406 times. The latter is a very high number for a mathematician. Since 1990 Klibanov is with University of North Carolina at Charlotte, USA. 
\end{IEEEbiography}
\begin{IEEEbiography}
[{\includegraphics[width=1in,height=1.25in,clip,keepaspectratio]{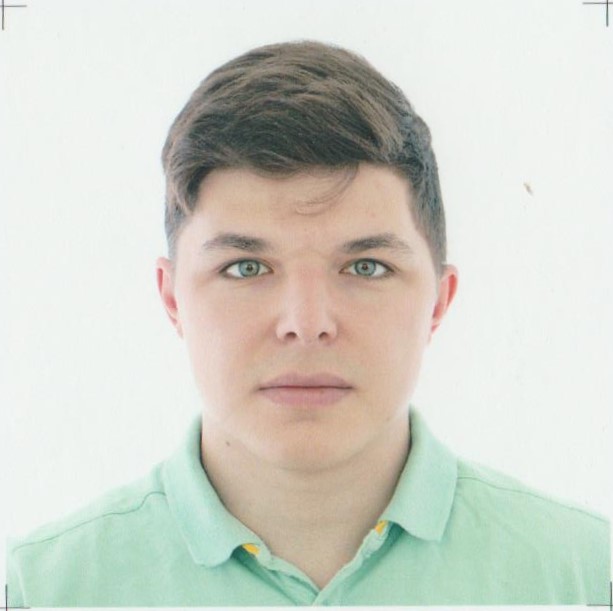}}]{Alexey V. Smirnov}
received the B.S. in physics from the Lomonosov Moscow State University, Moscow, Russia and the M.S. degree in mathematics from University of North Carolina at Charlotte, Charlotte, NC, USA in 2018 and 2020, respectively. He is currently working towards  a Ph.D. degree in applied mathematics at the University of Waterloo, Waterloo, ON, Canada. His current research interests include numerical methods for inverse problems with applications in medical imaging and remote sensing. 
\end{IEEEbiography}

\begin{IEEEbiography}
[{\includegraphics[width=1in,height=1.25in,clip,keepaspectratio]{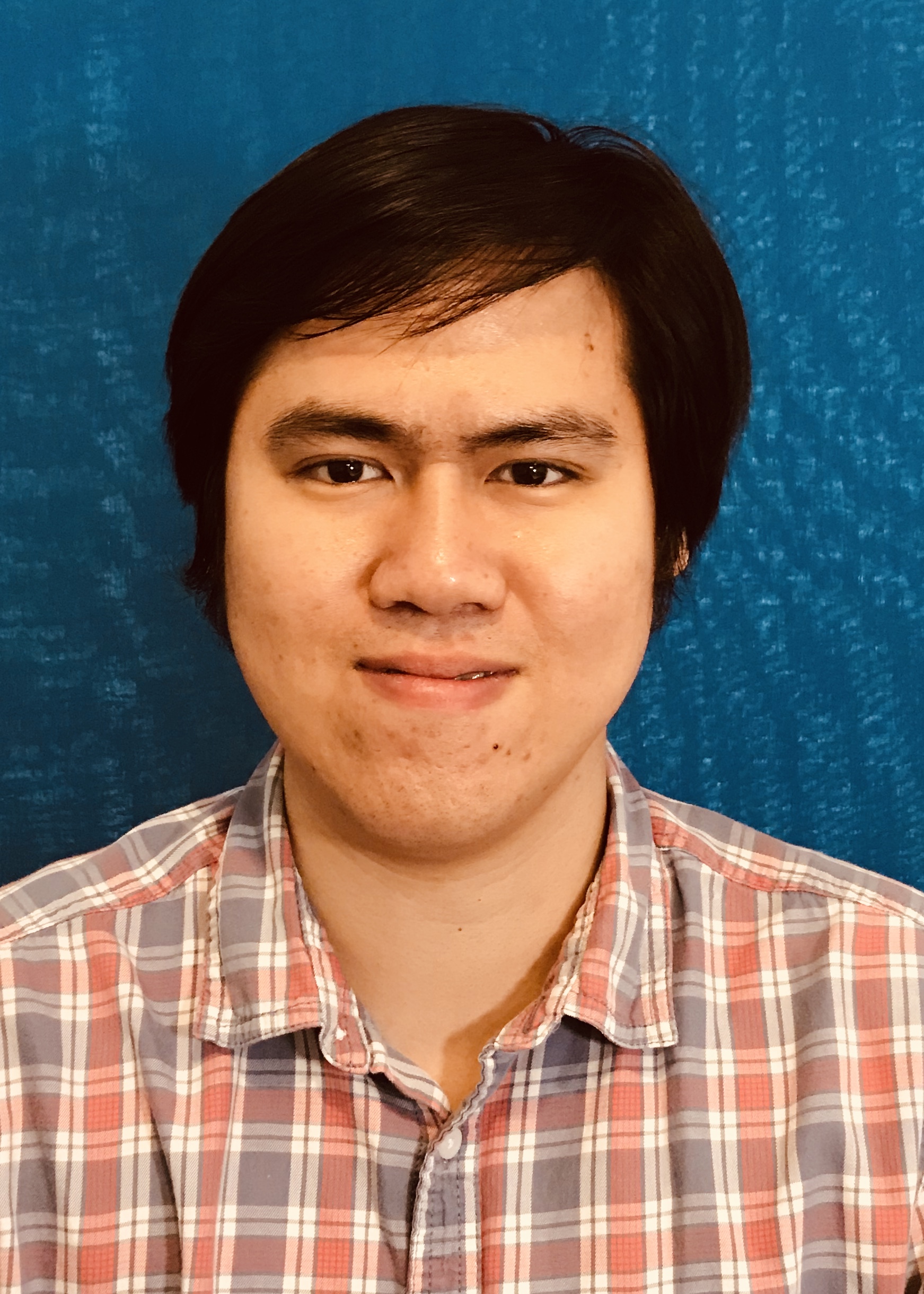}}]{Vo A. Khoa} received the B.Sc. degree (honor) in mathematics and computer Science from Ho Chi Minh City University of Science, Vietnam, in 2014. In 2018, he obtained the Ph.D. degree (cum laude) in mathematics from a joint international program between the International School for Advanced Studies (SISSA) and the Gran Sasso Science Institute in Italy. Since 2019, Dr. Khoa has served as a postdoctoral scholar at University of North Carolina at Charlotte, NC, USA. Between 2018 and 2019, he was a short-term postdoc at University of Goettingen, Germany, and then won a postdoctoral fellowship, hosted by Hasselt University, Belgium, from the Research Foundation - Flanders. His current research is specialized in inverse and ill-posed problems for partial differential equations, solving identification problems in physics and biology. Dr. Khoa is an active member of the inverse and ill-posed community. In 2017, he was among fifteen
outstanding reviewers of the prestigious journal Inverse Problems. Since 2019, he has been assigned to
volunteering as a reviewer for Mathematical Reviews of American Mathematical Society.
\end{IEEEbiography}
\begin{IEEEbiography}
[{\includegraphics[width=1in,height=1.25in,clip,keepaspectratio]{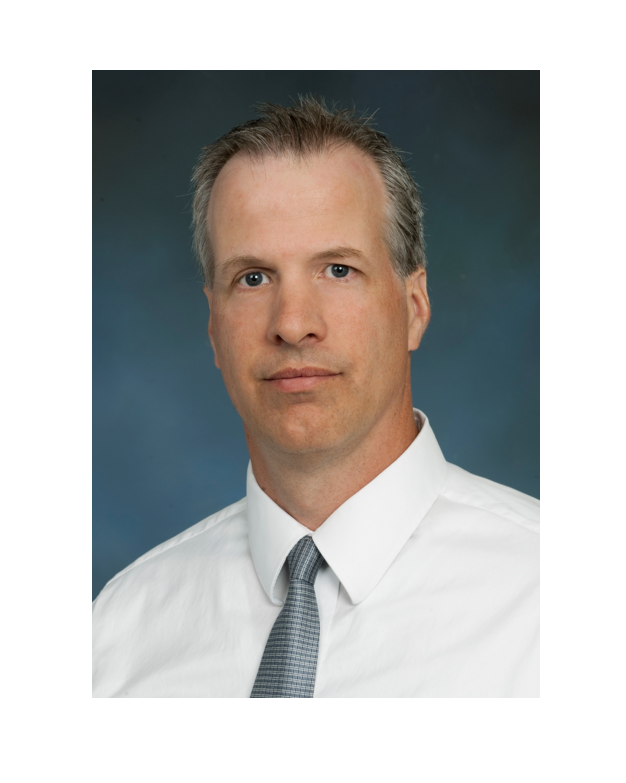}}]{Anders J. Sullivan}
received the B.S. and M.S. degrees in aerospace engineering from the Georgia Institute of Technology, Atlanta, and the Ph.D. degree from Polytechnic University, Brooklyn, NY, with a specialty in electromagnetics. He began his career with the U.S. Air Force Research Laboratory, Eglin Air Force Base, FL. Following this, he was a Postdoctoral Research Associate with the Electrical and Computer Engineering Department, Duke University, Durham, NC.  He is currently a branch chief at the Army Research Laboratory in Adelphi, MD. His main research interests include computational electromagnetics and signal processing associated with concealed target detection radar applications. 
\end{IEEEbiography}
\vfill
\begin{IEEEbiography}
[{\includegraphics[width=1in,height=1.25in,clip,keepaspectratio]{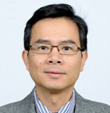}}]{Lam H. Nguyen}
received the B.S.E.E. degree from Virginia Polytechnic Institute, Blacksburg, VA, USA, the M.S.E.E. degree from George Washington University, Washington, DC, USA, and the M.S.C.S. degree from Johns Hopkins University, Baltimore, MD, USA, in 1984, 1991, and 1995, respectively. He started his career with General Electric Company, Portsmouth, VA, in 1984. He joined Harry Diamond Laboratory, Adelphi, MD (and its predecessor Army Research Laboratory) and has worked there from 1986 to the present. Currently, he is a Signal Processing Team Leader with the U.S. Army Research Laboratory, where he has primarily engaged in the research and development of several versions of ultra-wideband (UWB) radar since the early 1990s to the present. These radar systems have been used for proof-of-concept demonstrations in many concealed target detection programs. He has authored/coauthored approximately 100 conferences, journals, and technical publications. He has twelve patents in SAR system and signal processing. He has been a member of the SPIE Technical Committee on Radar Sensor Technology since 2009. He was the recipient of the U.S. Army Research and Development Achievement Awards in 2006,
2008, and 2010, the Army Research Laboratory Award for Science in 2010, and the U.S. Army Superior Civilian Performance Award in 2011.
\end{IEEEbiography}
\vfill


\end{document}